\documentclass[a4paper,12pt,twoside,leqno,final]{amsart}
\usepackage{amsmath}
\usepackage{amssymb}

\newtheorem{thm}{Theorem}[section]
\newtheorem{lem}[thm]{Lemma}

\newcommand{\C}{{\mathbb C}}
\newcommand{\D}{{\mathbb D}}

\newcommand{\T}{{\mathbb T}}

\newcommand{\N}{{\mathbb N}}

\newcommand{\bmo}{{\rm BMO}}
\newcommand{\bmoa}{{\rm BMOA}}

\newcommand{\eps}{\varepsilon}

\newcommand{\f}{\frac}
\newcommand{\ov}{\overline}
\newcommand{\al}{\alpha}
\newcommand{\be}{\beta}

\newcommand{\ga}{\gamma}

\newcommand{\de}{\delta}

\newcommand{\ze}{\zeta}
\renewcommand{\th}{\theta}

\newcommand{\ph}{\varphi}

\numberwithin{equation}{section}

\title[A free interpolation problem for a subspace of $H^\infty$]
{A free interpolation problem\\ 
for a subspace of $H^\infty$}
\author{Konstantin M. Dyakonov}
\address{Departament de Matem\`atiques i Inform\`atica, IMUB, BGSMath, Universitat de Barcelona, Gran Via 585, E-08007 Barcelona, Spain}
\address{ICREA, Pg. Llu\'is Companys 23, E-08010 Barcelona, Spain}
\email{konstantin.dyakonov@icrea.cat}
\thanks{Supported in part by grants MTM2014-51834-P and MTM2017-83499-P from El Ministerio de Econom\'ia y Competitividad (Spain), and by grant 2017-SGR-358 from AGAUR (Generalitat de Catalunya).}

\begin{document}
\begin{abstract}
Given an inner function $\th$, the associated star-invariant subspace $K^\infty_\th$ is formed by the functions $f\in H^\infty$ that annihilate (with respect to the usual pairing) the shift-invariant subspace $\th H^1$ of the Hardy space $H^1$. Assuming that $B$ is an interpolating Blaschke product with zeros $\{a_j\}$, we characterize the traces of functions from $K^\infty_B$ on the sequence $\{a_j\}$. The trace space that arises is, in general, non-ideal (i.e., the sequences $\{w_j\}$ belonging to it admit no nice description in terms of the size of $|w_j|$), but we do point out explicit -- and sharp -- size conditions on $|w_j|$ which make it possible to solve the interpolation problem $f(a_j)=w_j$ ($j=1,2,\dots$) with a function $f\in K^\infty_B$.
\end{abstract}

\maketitle

\section{Introduction and results}

Let $\D$ stand for the open unit disk (centered at $0$) in the complex plane, $\T$ for its boundary, and $m$ for the normalized arclength measure on $\T$. Recall, further, that the {\it Hardy space} $H^p=H^p(\D)$ with $0<p<\infty$ is defined as the set of all holomorphic functions $f$ on $\D$ that satisfy 
$$\sup_{0<r<1}\int_\T|f(r\ze)|^pdm(\ze)<\infty,$$
while $H^\infty=H^\infty(\D)$ denotes the space of bounded holomorphic functions. As usual, we identify elements of $H^p$ with their boundary functions (living almost everywhere on $\T$) and thus embed $H^p$ isometrically into $L^p=L^p(\T,m)$. The underlying theory can be found in \cite[Chapter II]{G}. 

\par Now suppose $\th$ is an {\it inner function} on $\D$; this means, by definition, that $\th\in H^\infty$ and $|\th|=1$ a.e. on $\T$. We then introduce the associated {\it star-invariant} (or {\it model}) {\it subspace} $K^p_\th$, this time with $1\le p\le\infty$, by putting 
\begin{equation}\label{eqn:defnkpth}
K^p_\th:=H^p\cap\th\,\ov{H^p_0}, 
\end{equation}
where $H^p_0:=\{f\in H^p:f(0)=0\}$ and the bar denotes complex conjugation. The term \lq\lq star-invariant" means invariant under the backward shift operator $$f\mapsto\f{f-f(0)}z,$$ 
and it is known (cf. \cite{DSS}) that the general form of a proper closed star-invariant subspace in $H^p$, with $1\le p<\infty$, is indeed given by \eqref{eqn:defnkpth} for some inner function $\th$. The alternative term \lq\lq model subspace" is due to the appearance of these subspaces in the Sz.-Nagy--Foia\c{s} operator model; see \cite{N}. It follows from the definition that each $K^p_\th$ carries a natural antilinear isometry (or involution) given by $f\mapsto\widetilde f$, where 
\begin{equation}\label{eqn:ftilde}
\widetilde f:=\ov z\ov f\th.
\end{equation}

\par When $p=2$, we can equivalently define $K^2_\th$ as the orthogonal complement of the shift-invariant subspace $\th H^2$ in $H^2$. Moreover, letting $P_\th$ denote the orthogonal projection from $H^2$ onto $K^2_\th$, one easily verifies that 
$$P_\th f=f-\th P_+(\ov\th f)=\th P_-(\ov\th f)$$
for $f\in H^2$, where $P_+$ and $P_-$ are the orthogonal projections from $L^2$ onto $H^2$ and $\ov{H^2_0}$, respectively. Now, we know from the M. Riesz theorem (see \cite[Chapter III]{G}) that $P_+$ and $P_-$ extend -- or restrict -- to every $L^p$ space with $1<p<\infty$ as bounded operators (called the {\it Riesz projections}), their respective ranges being $H^p$ and $\ov{H^p_0}$. It follows then that $P_\th$ admits a bounded extension -- or restriction -- to every $H^p$ with $1<p<\infty$, and projects the latter space onto $K^p_\th$ parallel to $\th H^p$. Accordingly, we arrive at the direct sum decomposition 
\begin{equation}\label{eqn:dirsumth}
H^p=K^p_\th\oplus\th H^p,\qquad1<p<\infty,
\end{equation}
with orthogonality for $p=2$. 

\par We shall make use of \eqref{eqn:dirsumth} in a special situation, which we now describe. Recall that, given a sequence $\{a_j\}=\{a_j\}_{j\in\N}$ of points in $\D$ with $\sum_j(1-|a_j|)<\infty$, the associated {\it Blaschke product} $B$ is defined by 
\begin{equation}\label{eqn:blapro}
B(z)=B_{\{a_j\}}(z):=\prod_jb_j(z),\quad\text{\rm where}\quad b_j(z):=\f{|a_j|}{a_j}\f{a_j-z}{1-\ov a_jz},
\end{equation}
with the convention that $|a_j|/a_j=-1$ if $a_j=0$. Then $B$ is an inner function that vanishes precisely at the $a_j$'s; see \cite[Chapter II]{G}. Recall also that a sequence $\{a_j\}$ in $\D$ is called an {\it interpolating sequence} if 
\begin{equation}\label{eqn:trcarl}
H^\infty\big|_{\{a_j\}}=\ell^\infty.
\end{equation}
(Here and below, given a function space $\mathcal X$ on $\D$, we denote by $\mathcal X\big|_{\{a_j\}}$ the set of those sequences $\{w_j\}$ in $\C$ for which the interpolation problem $f(a_j)=w_j$, $j\in\N$, has a solution $f\in\mathcal X$.) Carleson's celebrated theorem (see \cite{Carl} or \cite[Chapter VII]{G}) characterizes the interpolating sequences $\{a_j\}$ in terms of the corresponding Blaschke product \eqref{eqn:blapro} or rather its subproducts $B_j:=B/b_j$. Namely, it asserts that $\{a_j\}$ is interpolating if and only if 
\begin{equation}\label{eqn:carlcond}
\inf_j|B_j(a_j)|>0,
\end{equation}
a condition that can be further rephrased as 
\begin{equation}\label{eqn:intbla}
\inf_j|B'(a_j)|\,(1-|a_j|)>0.
\end{equation}
A Blaschke product $B=B_{\{a_j\}}$ satisfying \eqref{eqn:intbla} is said to be an {\it interpolating Blaschke product}. 

\par When $0<p<\infty$, we have a similar \lq\lq free interpolation" phenomenon in $H^p$. This time, \eqref{eqn:trcarl} gets replaced by 
\begin{equation}\label{eqn:trshsh}
H^p\big|_{\{a_j\}}=\left\{\{w_j\}:\,\sum_j|w_j|^p(1-|a_j|)<\infty\right\},
\end{equation}
and results of \cite{Kab, SS} tell us that this happens, for some or each $p\in(0,\infty)$, if and only if $\{a_j\}$ obeys the Carleson condition \eqref{eqn:carlcond}. Now, in the case $1<p<\infty$, we may apply \eqref{eqn:dirsumth} with $\th=B\left(=B_{\{a_j\}}\right)$, and restricting both sides to $\{a_j\}$ yields 
\begin{equation}\label{eqn:equaltraces}
H^p\big|_{\{a_j\}}=K^p_B\big|_{\{a_j\}}.
\end{equation}
Finally, we combine \eqref{eqn:trshsh} and \eqref{eqn:equaltraces} to conclude that 
\begin{equation}\label{eqn:tracekpb}
K^p_B\big|_{\{a_j\}}=\left\{\{w_j\}:\,\sum_j|w_j|^p(1-|a_j|)<\infty\right\},\qquad1<p<\infty,
\end{equation}
whenever $B$ is an interpolating Blaschke product with zeros $\{a_j\}$. 

\par For $p=\infty$, however, no such thing is true, since the endpoint version of \eqref{eqn:equaltraces} breaks down in general. A natural problem that arises is, therefore, to find out what happens to \eqref{eqn:trcarl} if we replace $H^\infty$ by its star-invariant subspace 
$$K^\infty_B:=H^\infty\cap B\ov{H^\infty_0},$$
always assuming that $B=B_{\{a_j\}}$ is an interpolating Blaschke product. 

\par It does happen sometimes that 
\begin{equation}\label{eqn:maxtrace}
K^\infty_B\big|_{\{a_j\}}=\ell^\infty, 
\end{equation}
but typically, and in \lq\lq most" cases, our trace space will be essentially smaller. As a matter of fact, \eqref{eqn:maxtrace} holds if and only if $\{a_j\}$ is an interpolating sequence that satisfies the so-called {\it (uniform) Frostman condition}: 
\begin{equation}\label{eqn:ufc}
\sup\left\{\sum_j\f{1-|a_j|}{|\ze-a_j|}:\,\ze\in\T\right\}<\infty.
\end{equation}
(This result is a fairly straightforward consequence of Hru\v s\v cev and Vinogradov's work in \cite{HV}; see also \cite[Section 3]{C2} for details.) In particular, \eqref{eqn:ufc} implies that the $a_j$'s may only approach the unit circle in a suitably tangential manner. 

\par On the other hand, it was shown by Vinogradov in \cite{V} that whenever $\{a_j\}$ is an interpolating sequence, one has 
$$K^\infty_{B^2}\big|_{\{a_j\}}=\ell^\infty,$$
where again $B=B_{\{a_j\}}$. Note, however, that $K^\infty_{B^2}$ is substantially larger than $K^\infty_B$. 

\par We now describe the trace space $K^\infty_B\big|_{\{a_j\}}$ in the general case.

\begin{thm}\label{thm:interpolthm} Suppose that $\{w_j\}\in\ell^\infty$ and $B$ is an interpolating Blaschke product with zeros $\{a_j\}$. In order that 
\begin{equation}\label{eqn:wbeltrace}
\{w_j\}\in K^\infty_B\big|_{\{a_j\}},
\end{equation}
it is necessary and sufficient that 
\begin{equation}\label{eqn:crucond}
\sup_k\left|\sum_j\f{w_j}{B'(a_j)\cdot(1-a_j\ov a_k)}\right|<\infty.
\end{equation}
\end{thm}

An equivalent formulation is as follows. 

\begin{thm}\label{thm:equivform} Suppose $h\in H^\infty$ and $B$ is an interpolating Blaschke product with zeros $\{a_j\}$. Then $h\in K^\infty_B+BH^\infty$ if and only if 
$$\sup_k\left|\sum_j\f{h(a_j)}{B'(a_j)\cdot(1-a_j\ov a_k)}\right|<\infty.$$ 
\end{thm}

To deduce Theorem \ref{thm:equivform} from Theorem \ref{thm:interpolthm} and vice versa, it suffices to observe that a given $H^\infty$ function will be in $K^\infty_B+BH^\infty$ if and only if its values at the $a_j$'s can be interpolated by those of a function in $K^\infty_B$. 

\par Before stating our next result, we need to recall the notion of an ideal sequence space. A vector space $S$ consisting of sequences (of complex numbers) is said to be {\it ideal} if, whenever $\{x_j\}\in S$ and $\{y_j\}$ is a sequence with $|y_j|\le|x_j|$ for all $j$, we necessarily have $\{y_j\}\in S$. Roughly speaking, this property tells us that the elements $\{x_j\}$ of $S$ are somehow describable by means of a \lq\lq size condition" on $|x_j|$. 

\par The trace space in \eqref{eqn:tracekpb} is obviously ideal, but its endpoint version $K^\infty_B\big|_{\{a_j\}}$ can no longer be expected to have this nice feature. Of course, the latter space {\it will} be ideal for the \lq\lq few" interpolating sequences $\{a_j\}$ that obey the Frostman condition \eqref{eqn:ufc}, in which case we have \eqref{eqn:maxtrace}, but other choices of $\{a_j\}$ make things different. The difference becomes especially dramatic in the \lq\lq anti-Frostman" situation where the $a_j$'s are taken to lie on a radius. In our next theorem, we furnish a {\it universal} ideal sequence space, namely $\ell^1$, that is contained in every trace space $K^\infty_B\big|_{\{a_j\}}$, and we show (by examining the radial case) that no larger ideal space would do in general. 

\begin{thm}\label{thm:elloneideal} Suppose $B$ is an interpolating Blaschke product with zeros $\{a_j\}$. 
\par{\rm (a)} We have 
\begin{equation}\label{eqn:ellonetrace}
\ell^1\subset K^\infty_B\big|_{\{a_j\}}.
\end{equation}
\par{\rm (b)} If $0\le a_1<a_2<\dots<1$, and if $X$ is an ideal sequence space with 
\begin{equation}\label{eqn:xtrace}
X\subset K^\infty_B\big|_{\{a_j\}},
\end{equation}
then $X\subset\ell^1$.
\end{thm}

\par At the same time, it is worth mentioning that the trace space $K^\infty_B\big|_{\{a_j\}}$ always contains non-$\ell^1$ sequences (assuming, as we do, that $B=B_{\{a_j\}}$ is an {\it infinite} Blaschke product). An example is provided by the constant sequence consisting of $1$'s, since these are the values of the function $1-\ov{B(0)}B\,\left(\in K^\infty_B\right)$ at the $a_j$'s. 

\par We shall be also concerned with uniform smallness conditions on the values $w_j$ that guarantee \eqref{eqn:wbeltrace}, once the $a_j$'s are given. To be more precise, we fix a (reasonable) positive function $\ph$ on the interval $(0,1]$ and write $X_\ph(\{a_j\})$ for the set of all sequences $\{w_j\}\in\ell^\infty$ that satisfy   
\begin{equation}\label{eqn:sizewj}
\sup_j\f{|w_j|}{\ph(1-|a_j|)}<\infty.
\end{equation}
Our aim is then to determine whether 
\begin{equation}\label{eqn:xphintrace}
X_\ph(\{a_j\})\subset K^\infty_B\big|_{\{a_j\}} 
\end{equation}
with $B=B_{\{a_j\}}$, for every interpolating sequence ${\{a_j\}}$. This will be settled by Theorem \ref{thm:phiideal} below, but first we have to describe the class of $\ph$'s we have in mind. 

\par It will be assumed that $\ph:(0,1]\to(0,\infty)$ is a nondecreasing continuous function for which $t\mapsto\ph(t)/t$ is nonincreasing; a function $\ph$ with these properties will be called a {\it majorant}. Also, following the terminology of \cite{J}, we say that $\ph$ is {\it of upper type less than $1$} if there are constants $C>0$ and $\ga\in[0,1)$ such that 
\begin{equation}\label{eqn:uppertype}
\ph(st)\le Cs^\ga\ph(t)
\end{equation}
whenever $s\ge1$ and $0<t\le1/s$. It should be noted that for $\ga=1$, \eqref{eqn:uppertype} is automatic (with $C=1$) for every majorant $\ph$. 

\begin{thm}\label{thm:phiideal} {\rm (i)} If $\ph$ is a majorant of upper type less than $1$ satisfying 
\begin{equation}\label{eqn:intconv}
\int_0^1\f{\ph(t)}t\,dt<\infty,
\end{equation}
then, for every interpolating Blaschke product $B=B_{\{a_j\}}$, we have \eqref{eqn:xphintrace}. 
\par{\rm (ii)} Conversely, if $\ph$ is a majorant with the property that \eqref{eqn:xphintrace} is valid for each interpolating Blaschke product $B=B_{\{a_j\}}$, then \eqref{eqn:intconv} holds true. 
\end{thm}

In fact, a glance at the proof of part (ii) will reveal that it is enough to assume \eqref{eqn:xphintrace} for a {\it single} interpolating sequence $\{a_j\}$, namely, for $a_j=1-2^{-j}$; this alone will imply \eqref{eqn:intconv}. 
\par As examples of majorants $\ph$ that are of upper type less than $1$ and satisfy \eqref{eqn:intconv}, one may consider 
$$\ph_1(t)=t^\al,\quad\ph_2(t)=\left(\log\f2t\right)^{-1-\eps},\quad \ph_3(t)=\left(\log\f3t\right)^{-1}\left(\log\log\f3t\right)^{-1-\eps}$$
(with $0<\al<1$ and $\eps>0$), and so on. 

\par While we are only concerned with the traces of functions from $K^\infty_B$ on $\{a_j\}$, which is the zero sequence of $B$, an obvious generalization would consist in restricting our functions (or those from $K^\infty_\th$, with $\th$ inner) to an arbitrary interpolating sequence in $\D$. In this generality, however, even the case of $K^p_\th$ with $1<p<\infty$ (or with $p=2$) cannot be viewed as completely understood. At least, the existing results (see \cite{AH, DPAMS, HNP} for some of these) appear to be far less clear-cut than in the current setting. By contrast, the difficulty we have to face here is entirely due to the endpoint position of $K^\infty_B$ within the $K^p_B$ scale, the only enemy being the failure of the M. Riesz projection theorem. 

\par The other endpoint case, where $p=1$, is no less wicked and we briefly discuss it here. For an interpolating sequence $\{a_k\}$, the values $w_k=f(a_k)$ of a function $f\in H^1$ satisfy 
\begin{equation}\label{eqn:carlhone}
\sum_k|w_k|\,(1-|a_k|)<\infty
\end{equation}
and, in virtue of the Shapiro--Shields theorem \cite{SS}, this property characterizes the sequences $\{w_k\}$ in $H^1\big|_{\{a_k\}}$. The latter set will, however, be strictly larger than $K^1_B\big|_{\{a_k\}}$ with $B=B_{\{a_k\}}$, unless we are dealing with finitely many $a_k$'s. (Equivalently, we never have $K^1_B\oplus BH^1=H^1$ except when $B$ is a finite Blaschke product. This can be deduced from \cite[Theorem 3.8]{Steg}.) Thus, in the nontrivial cases, \eqref{eqn:carlhone} is necessary but not sufficient for $\{w_k\}$ to be the trace of a function from $K^1_B$ on $\{a_k\}$. Now, an adaptation of our current method leads to another necessary condition involving the numbers 
$$\widetilde w_k:=\sum_j\f{w_j}{B'(a_j)\cdot(1-a_j\ov a_k)},\qquad k\in\N.$$
Namely, \eqref{eqn:carlhone} must also hold with $\widetilde w_k$ in place of $w_k$. It would be interesting to determine whether the two conditions together are actually sufficient for $\{w_k\}$ to be in $K^1_B\big|_{\{a_k\}}$. A more detailed discussion and further questions can be found in \cite{DIEOT}. 

\par In the next section, we collect a few auxiliary facts to lean upon. The remaining sections contain the proofs of our results.

\section{Preliminaries}

Given an inner function $\th$, we write
$$K_{*\th}:=K^2_\th\cap\bmo,$$
where $\bmo=\bmo(\T)$ is the space of functions of bounded mean oscillation on $\T$; see \cite[Chapter VI]{G}. The following representation formula is borrowed from \cite[Lemma 3.1]{C1}. 

\begin{lem}\label{lem:genformkbmo} Given an interpolating Blaschke product $B$ with zeros $\{a_j\}$, the general form of a function $g\in K_{*B}$ is 
\begin{equation}\label{eqn:kbmoseries}
g(z)=\sum_jc_j\f{1-|a_j|}{1-\ov a_jz}
\end{equation}
with $\{c_j\}\in\ell^\infty$.
\end{lem}

The series in \eqref{eqn:kbmoseries} is understood to converge in the weak-star topology of $\bmoa:=\bmo\cap H^2$, viewed as the dual of $H^1$. It also converges in $H^2$ (cf. \cite{HNP}), and hence on compact subsets of $\D$. 
\par We now cite another result of Cohn, namely \cite[Corollary 3.2]{C2}, as our next lemma. 

\begin{lem}\label{lem:kbmobdd} If $B$ is an interpolating Blaschke product with zeros $\{a_j\}$, and if $g\in K_{*B}$ is a function satisfying $\{g(a_j)\}\in\ell^\infty$, then $g\in H^\infty$.
\end{lem}

Further, we need to recall the definition of the space $\bmo_\ph=\bmo_\ph(\T)$ associated with a majorant $\ph$. A function $f\in L^1(\T,m)$ is said to be in $\bmo_\ph$ if 
$$\sup_I\f1{m(I)\cdot\ph(m(I))}\int_I|f-f_I|\,dm<\infty,$$
where $f_I:=m(I)^{-1}\int_If\,dm$, the supremum being taken over the open arcs $I\subset\T$. The classical $\bmo$ corresponds to the constant majorant $\ph\equiv1$. 
\par The following fact (and its converse) can be found in \cite{Sp}. 

\begin{lem}\label{lem:bmocont} If $\ph$ is a majorant satisfying \eqref{eqn:intconv}, then $\bmo_\ph\subset C(\T)$.
\end{lem}

Our last lemma is essentially contained in \cite{DJAM98}. 

\begin{lem}\label{lem:hankelbmo} Suppose that $B$ is an interpolating Blaschke product with zeros $\{a_j\}$, $\ph$ is a majorant of upper type less than $1$, and $f\in H^2$. If 
\begin{equation}\label{eqn:sizefataj}
\sup_j\f{|f(a_j)|}{\ph(1-|a_j|)}<\infty,
\end{equation}
then $P_-(f\ov B)\in\bmo_\ph$.
\end{lem}

Precisely speaking, this result was incorporated into the proof of Theorem 2.3 in \cite{DJAM98}. The theorem asserted (among other things) that \eqref{eqn:sizefataj} is necessary and sufficient in order that $f\ov B\in\bmo_\ph$, provided that $f\in\bmoa_\ph:=\bmo_\ph\cap H^2$. The sufficiency was then established by splitting $f\ov B$ as 
$$P_+(f\ov B)+P_-(f\ov B)$$ 
and verifying that both terms are in $\bmo_\ph$. In particular, the second term, $P_-(f\ov B)$, was shown to be in $\bmo_\ph$ by means of a duality argument (based on a result from \cite{J}) which actually works for any $f\in H^2$, not just for $f\in\bmoa_\ph$; see \cite[p.\,97]{DJAM98} for details. 
\par When $\ph(t)=t^\al$ with some $\al\in(0,1)$, $\bmo_\ph$ reduces to the usual Lipschitz space of order $\al$, and in this special case Lemma \ref{lem:hankelbmo} was previously established in \cite[Section 4]{DSpb93}; the case of higher order Lipschitz--Zygmund spaces was treated there as well. We also refer to \cite{DIUMJ94, DAiM17} for related results.

\section{Proof of Theorem \ref{thm:interpolthm}}

Suppose \eqref{eqn:wbeltrace} holds, so that there exists a function $f\in K^\infty_B$ satisfying
\begin{equation}\label{eqn:valuefaj}
f(a_j)=w_j,\qquad j\in\N.
\end{equation}
To deduce \eqref{eqn:crucond}, we first define 
\begin{equation}\label{eqn:defgammaj}
\ga_j:=-\f{a_j}{|a_j|}\f{w_j}{B_j(a_j)},\qquad j\in\N
\end{equation}
(where $B_j$ is the Blaschke product with zeros $\{a_k:k\ne j\}$), and consider the function 
\begin{equation}\label{eqn:functg}
g(z):=\sum_j\ov\ga_j\f{1-|a_j|^2}{1-\ov a_jz}.
\end{equation}
Observe that $\{\ga_j\}\in\ell^\infty$, because $\inf_j|B_j(a_j)|>0$, and so $g\in K_{*B}$ by virtue of Lemma \ref{lem:genformkbmo}; the latter is being applied with $c_j=\ov\ga_j(1+|a_j|)$. 
\par Recalling the notation \eqref{eqn:ftilde}, which is henceforth used with $B$ in place of $\th$, we have then (a.e. on $\T$)
\begin{equation*}
\begin{aligned}
\widetilde g(z)&:=\ov z\ov{g(z)}B(z)=\ov z\sum_jB_j(z)b_j(z)\ga_j\f{1-|a_j|^2}{1-a_j\ov z}\\
&=\ov z\sum_jB_j(z)\f{|a_j|}{a_j}\f{a_j-z}{1-\ov a_jz}
\ga_j\f{1-|a_j|^2}{1-a_j\ov z}\\
&=-\sum_jB_j(z)\f{|a_j|}{a_j}\ga_j\f{1-|a_j|^2}{1-\ov a_jz}\\
&=\sum_jB_j(z)\f{w_j}{B_j(a_j)}\f{1-|a_j|^2}{1-\ov a_jz}.
\end{aligned}
\end{equation*}
The resulting identity 
$$\widetilde g(z)=\sum_jB_j(z)\f{w_j}{B_j(a_j)}\f{1-|a_j|^2}{1-\ov a_jz}$$
actually holds for all $z\in\D$, and computing both sides at $a_k$ gives 
\begin{equation}\label{eqn:gtildeatak}
\widetilde g(a_k)=w_k,\qquad k\in\N
\end{equation}
(just note that $B_j(a_k)=0$ whenever $j\ne k$). 

\par Comparing \eqref{eqn:gtildeatak} and \eqref{eqn:valuefaj}, we deduce that $\widetilde g=f$; indeed, the difference $\widetilde g-f$ belongs to both $K^2_B$ and $BH^2$, and is therefore null. Because $f$ is actually in $K^\infty_B$, so is $g(=\widetilde f)$, and this obviously implies that 
\begin{equation}\label{eqn:gatakbdd} 
\{g(a_k)\}\in\ell^\infty.
\end{equation}
Equivalently, the values $\ov{g(a_k)}$ ($k\in\N$) form a bounded sequence, i.e., 
\begin{equation}\label{eqn:bargatak}
\sup_k\left|\sum_j\ga_j\f{1-|a_j|^2}{1-a_j\ov a_k}\right|<\infty.
\end{equation}
Finally, we combine \eqref{eqn:defgammaj} with the elementary formula 
$$B'(a_j)=B_j(a_j)\cdot b'_j(a_j)=-\f{|a_j|}{a_j}\f{B_j(a_j)}{1-|a_j|^2}$$
to get 
\begin{equation}\label{eqn:uside}
\ga_j\cdot(1-|a_j|^2)=\f{w_j}{B'(a_j)},
\end{equation}
and substituting this into \eqref{eqn:bargatak} yields \eqref{eqn:crucond}. 

\par Conversely, assume that \eqref{eqn:crucond} holds. Further, let $f\in K^2_B$ be a function satisfying \eqref{eqn:valuefaj}. (To find such an $f$, it suffices to solve the interpolation problem with an $H^2$ function and then project it orthogonally onto $K^2_B$.) Defining the numbers $\ga_j$ and the function $g$ by \eqref{eqn:defgammaj} and \eqref{eqn:functg}, respectively, we then infer -- exactly as before -- that $g$ lies in $K_{*B}$ and obeys \eqref{eqn:gtildeatak}. The latter in turn implies that $\widetilde g=f$, as above. 

\par On the other hand, we may again use the identity \eqref{eqn:uside} to rewrite \eqref{eqn:crucond} as \eqref{eqn:bargatak}, or equivalently, as \eqref{eqn:gatakbdd}. This done, we invoke Lemma \ref{lem:kbmobdd} to conclude that $g$ is in $H^\infty$, and hence in $K^\infty_B$. The function $f(=\widetilde g)$ therefore belongs to $K^\infty_B$ as well, and recalling \eqref{eqn:valuefaj} we finally arrive at \eqref{eqn:wbeltrace}. 

\section{Proofs of Theorems \ref{thm:elloneideal} and \ref{thm:phiideal}}

\noindent{\it Proof of Theorem \ref{thm:elloneideal}.} (a) For all $j$ and $k$ in $\N$, we have 
$$|B'(a_j)|\cdot|1-a_j\ov a_k|\ge|B'(a_j)|\cdot(1-|a_j|)\ge\de>0,$$ 
where $\de$ is the infimum in \eqref{eqn:intbla}. Therefore, whenever $\{w_j\}\in\ell^1$, 
$$\sup_k\left|\sum_j\f{w_j}{B'(a_j)\cdot(1-a_j\ov a_k)}\right|
\le\f1\de\sum_j|w_j|<\infty,$$
and \eqref{eqn:crucond} holds true. In view of Theorem \ref{thm:interpolthm}, the inclusion \eqref{eqn:ellonetrace} is thereby verified. 

\smallskip (b) Assume, under the current hypotheses on $\{a_j\}$ and $X$, that we can find a sequence $\{\be_j\}\in X\setminus\ell^1$. Put 
$$w_j:=B'(a_j)\cdot(1-a^2_j)\cdot|\be_j|,\qquad j\in\N.$$ 
Because $B$ is a unit-norm $H^\infty$ function, we have 
$$|B'(a_j)|\cdot(1-a^2_j)\le1$$
(we are also using the fact that $0\le a_j<1$ for all $j$), and so 
$$|w_j|\le|\be_j|,\qquad j\in\N.$$
Since $X$ is an ideal sequence space containing $\{\be_j\}$, it follows that $\{w_j\}\in X$. Recalling \eqref{eqn:xtrace}, we readily arrive at \eqref{eqn:wbeltrace}, which we further rephrase (using Theorem \ref{thm:interpolthm}), as \eqref{eqn:crucond}. Thus, the sums 
$$S_k:=\sum_j\f{w_j}{B'(a_j)\cdot(1-a_ja_k)}$$
(whose terms are all real and nonnegative) must satisfy 
\begin{equation}\label{eqn:crucondspcase}
\sup_kS_k<\infty.
\end{equation}

\par On the other hand, for any fixed $k\in\N$ and any $j\le k$, we have $a_k\ge a_j$, whence 
$$S_k\ge\sum_{j=1}^k\f{w_j}{B'(a_j)\cdot(1-a_ja_k)}
\ge\sum_{j=1}^k\f{w_j}{B'(a_j)\cdot(1-a^2_j)}=\sum_{j=1}^k|\be_j|.$$ 
Now, since $\sum_{j=1}^k|\be_j|\to\infty$ as $k\to\infty$, we conclude that $\sup_kS_k=\infty$, which contradicts \eqref{eqn:crucondspcase}. The contradiction means that the difference $X\setminus\ell^1$ is actually empty, so $X\subset\ell^1$ as required. \qed

\medskip\noindent{\it Proof of Theorem \ref{thm:phiideal}.} (i) Let $\{w_j\}$ be a sequence in $X_\ph(\{a_j\})$, so that \eqref{eqn:sizewj} holds, and let $f\in K^2_B$ be a solution to the interpolation problem \eqref{eqn:valuefaj}. Rewriting \eqref{eqn:sizewj} as \eqref{eqn:sizefataj}, while noting that 
$$f\ov B=P_-(f\ov B),$$
we deduce from Lemma \ref{lem:hankelbmo} that $f\ov B\in\bmo_\ph$. On the other hand, condition \eqref{eqn:intconv} guarantees, by virtue of Lemma \ref{lem:bmocont}, that every function in $\bmo_\ph$ is continuous (and hence bounded) on $\T$. It follows that $f\ov B$ is in $L^\infty$, and so is $f$. Consequently, $f$ actually belongs to $K^2_B\cap L^\infty=K^\infty_B$, and we arrive at \eqref{eqn:wbeltrace}. The inclusion \eqref{eqn:xphintrace} is thus established. 

\smallskip (ii) We put $a_j=1-2^{-j}$ ($j=1,2,\dots$) and exploit \eqref{eqn:xphintrace} in this special case only. The sequence space $X_\ph(\{a_j\})$ is obviously ideal, so we infer from Theorem \ref{thm:elloneideal}, part (b), that $X_\ph(\{a_j\})\subset\ell^1$; and since the sequence $\{\ph(1-a_j)\}=\{\ph(2^{-j})\}$ belongs to $X_\ph(\{a_j\})$, it follows in particular that 
\begin{equation}\label{eqn:phiellone}
\sum_{j=1}^\infty\ph(2^{-j})<\infty.
\end{equation}
It remains to observe that \eqref{eqn:phiellone} is equivalent to \eqref{eqn:intconv}. To see why, one only needs to break up the integration interval $(0,1]$ as $\bigcup_{j=1}^\infty I_j$, where $I_j:=(2^{-j},2^{-j+1}]$, and notice that 
$$\ph(2^{-j})\le\ph(t)\le2\ph(2^{-j})$$
for $t\in I_j$. The proof is complete. \qed


\medskip
 
\begin{thebibliography}{12} 

\bibitem{AH} E. Amar and A. Hartmann, {\it Uniform minimality, unconditionality and interpolation in backward shift invariant subspaces}, Ann. Inst. Fourier (Grenoble) \textbf{60} (2010), 1871--1903.

\bibitem{Carl} L. Carleson, {\it An interpolation problem for bounded analytic functions}, Amer. J. Math. \textbf{80} (1958), 921--930. 

\bibitem{C1} W. S. Cohn, {\it Radial limits and star invariant subspaces of bounded mean oscillation}, Amer. J. Math. \textbf{108} (1986), 719--749. 

\bibitem{C2} W. S. Cohn, {\it A maximum principle for star invariant subspaces}, Houston J. Math. \textbf{14} (1988), 23--37. 

\bibitem{DSS} R. G. Douglas, H. S. Shapiro, and A. L. Shields, {\it Cyclic vectors and invariant subspaces for the backward shift operator}, Ann. Inst. Fourier (Grenoble) \textbf{20} (1970), 37--76. 

\bibitem{DPAMS} K. M. Dyakonov, {\it Interpolating functions of minimal norm, star-invariant subspaces and kernels of Toeplitz operators}, Proc. Amer. Math. Soc. \textbf{116} (1992), 1007--1013. 

\bibitem{DSpb93} K. M. Dyakonov, {\it Smooth functions and coinvariant subspaces of the shift operator}, Algebra i Analiz \textbf{4} (1992), no. 5, 117--147; translation in St. Petersburg Math. J. \textbf{4} (1993), 933--959. 

\bibitem{DIUMJ94} K. M. Dyakonov, {\it Smooth functions in the range of a Hankel operator}, Indiana Univ. Math. J. \textbf{43} (1994), 805--838. 

\bibitem{DJAM98} K. M. Dyakonov, {\it Multiplicative structure in weighted $\bmoa$ spaces}, J. Anal. Math. \textbf{75} (1998), 85--104.

\bibitem{DAiM17} K. M. Dyakonov, {\it Factorization and non-factorization theorems for pseudocontinuable functions}, Adv. Math. \textbf{320} (2017), 630--651. 

\bibitem{DIEOT} K. M. Dyakonov, {\it Interpolating by functions from model subspaces in $H^1$}, Integral Equ. Oper. Theory, to appear; see also arXiv:1802.06433

\bibitem{G} J. B. Garnett, {\it Bounded analytic functions, Revised first edition}, Springer, New York, 2007. 

\bibitem{HNP} S. V. Hru\v s\v c\"ev, N. K. Nikol'ski\u\i, and B. S. Pavlov, {\it Unconditional bases of exponentials and of reproducing kernels}, in: Complex analysis and spectral theory (Leningrad, 1979/1980), pp. 214--335, Lecture Notes in Math., 864, Springer, Berlin and New York, 1981. 

\bibitem{HV} S. V. Hru\v s\v cev and S. A. Vinogradov, {\it Inner functions and multipliers of Cauchy type integrals}, Ark. Mat. \textbf{19} (1981), 23--42. 

\bibitem{J} S. Janson, {\it Generalizations of Lipschitz spaces and an application to Hardy spaces and bounded mean oscillation}, Duke Math. J. \textbf{47} (1980), 959--982.

\bibitem{Kab} V. Kabaila, {\it Interpolation sequences for $H_p$ classes in the case $p<1$}, Litovsk. Mat. Sb. \textbf{3} (1963), no. 1, 141--147. (Russian)

\bibitem{N} N. K. Nikolski, {\it Operators, Functions, and Systems: An Easy Reading, Vol. 2: Model operators and systems}, Mathematical Surveys and Monographs, 93, Amer. Math. Soc., Providence, RI, 2002. 

\bibitem{SS} H. S. Shapiro and A. L. Shields, {\it On some interpolation problems for analytic functions}, Amer. J. Math. \textbf{83} (1961), 513--532. 

\bibitem{Sp} S. Spanne, {\it Some function spaces defined using the mean oscillation over cubes}, Ann. Scuola Norm. Sup. Pisa Sci. Fis. Mat. (3) \textbf{19} (1965), 593--608. 

\bibitem{Steg} D. A. Stegenga, {\it Some bounded Toeplitz operators on $H^1$ and applications of the duality between $H^1$ and the functions of bounded mean oscillation}, Amer. J. Math. \textbf{98} (1976), 573--589. 

\bibitem{V} S. A. Vinogradov, {\it Some remarks on free interpolation by bounded and slowly growing analytic functions}, Zap. Nauchn. Sem. Leningrad. Otdel. Mat. Inst. Steklov. (LOMI) \textbf{126} (1983), 35--46. 

\end{thebibliography}
\end{document}